\documentclass[letterpaper,10
pt]{article}

\usepackage{fullpage}
\usepackage{enumerate}
\usepackage{authblk}
\usepackage[colorinlistoftodos]{todonotes}
\usepackage{verbatim}

\usepackage{section, amsthm, textcase, setspace, amssymb, lineno, amsmath, amssymb, amsfonts, latexsym, fancyhdr, longtable, ulem}
\usepackage{epsfig, graphicx, pstricks,pst-grad,pst-text,tikz, tkz-berge,colortbl}

\usepackage{float}
\usepackage{mathrsfs}

\usepackage{fullpage}
\usepackage{enumerate}
\usepackage{authblk}
\usepackage[colorinlistoftodos]{todonotes}
\usepackage{verbatim}
\usepackage{amsmath}

\usepackage{section, amsthm, textcase, setspace, amssymb, lineno, amsmath, amssymb, amsfonts, latexsym, fancyhdr, longtable, ulem}
\usepackage{epsfig, graphicx, pstricks,pst-grad,pst-text,tikz, tkz-berge,colortbl}

\usepackage{float}
\usepackage{mathrsfs}

\usepackage{mathtools}
\newtheorem{theorem}{Theorem}[section]

\DeclarePairedDelimiter{\ceil}{\lceil}{\rceil}
\DeclarePairedDelimiter{\floor}{\lfloor}{\rfloor}

\usetikzlibrary{graphs,graphs.standard,positioning,decorations.pathreplacing}

\renewcommand{\Affilfont}{\itshape\small}

\newcommand{\R}{\mathbb R}
\newcommand{\C}{\mathbb C}
\newcommand{\N}{\mathbb N}
\newcommand{\Q}{\mathbb Q}
\newcommand{\Z}{\mathbb Z}

\newtheorem{corollary}[theorem]{Corollary}
\newtheorem{lemma}[theorem]{Lemma}
\newtheorem{remark}[theorem]{Remark}

\begin{document}


\author{Vincent E. Coll, Jr.$^*$ and Lee B. Whitt$^{**}$}

\title{The classification of flat Riemannian metrics on the plane}

\maketitle

\noindent
\textit{$^*$Department of Mathematics, Lehigh University, Bethlehem, PA: vec208@lehigh.edu \\ $^{**}$Northrop Grumman - Distinguished Technical Fellow \rm{(}$retired$\rm{)}, $San~Diego, CA:lee.barlow.whitt@hotmail.com$}

\begin{abstract}
\noindent We classify all smooth flat Riemannian metrics on the two-dimensional plane. 
In the complete case, it is well-known that these metrics are isometric to the Euclidean metric.  In the incomplete case, there is an abundance of naturally-arising, non-isometric metrics that are relevant and useful. 
Remarkably, the study and classification of all flat Riemannian metrics on the plane -- as a subject -- is new to the literature. Much of our research focuses on conformal metrics of the form $e^{2\varphi}g_0$, where $\varphi: \R^2\rightarrow \R$ is a harmonic function and $g_0$ is the standard Euclidean metric on $\R^2$.  We find that all such metrics, which we call ``harmonic'', arise from Riemann surfaces. 

\end{abstract}



\section{Motivation for Incomplete Riemannian Metrics}

There is a broad consensus among mathematicians that incomplete Riemannian metrics are uninteresting. Reasons include the dominance of completeness as a condition for most of the important results in Riemannian geometry and the ease with which incomplete metrics can be 
contrived -- for example, remove a point from a Riemannian manifold.  

In this paper, we intend to challenge this consensus by offering a large family of naturally-occurring incomplete Riemannian metrics whose underlying geometries are non-trivial.  Surprisingly, this family consists of flat metrics on $\R^2$ and it is noteworthy that their classification fills a gap in the literature.


To begin, consider conformal metrics of the form $g = e^{2\varphi}g_0$, where $\varphi: \R^2 \rightarrow \R$ is harmonic and $g_0$ is the Euclidean metric on $\R^2$.  The Gaussian curvature of $g$ is easily calculated to be $K = -e^{2\varphi}\Delta_0\varphi$, where $\Delta_0$ is the standard Laplacian for $g_0$.  If follows then that $g$ is flat precisely when $\varphi$ is harmonic.  When  $\varphi$ is a non-constant harmonic function, we call $e^{2\varphi}g_0$ a \textit{harmonic metric} and the resulting flat plane a \textit{harmonic plane}.

Here are key facts about harmonic metrics and harmonic planes:

\begin{itemize}
\item Every harmonic metric on the plane is incomplete (see Lemma~\ref{incomplete}),

\item No harmonic plane is isometric to a proper subset of the Euclidean plane (see Lemma~\ref{another}).  More generally, no harmonic plane is isometric to a proper subset of another harmonic plane (see Lemma~\ref{another2}). 

\end{itemize}

One of the key reasons motivating the study of harmonic metrics is the direct correspondence between harmonic functions and harmonic metrics.  Preliminary investigations here reveal unexpected relationships between the analytic behavior of harmonic functions and the geometry e.g., the geodesics, of the associated harmonic metrics  
(see Section 7.4 and the examples in Section 8).

Here are two more reasons to motivate interest: 

\begin{enumerate}

\item Every harmonic plane is a Riemann surface.  The easiest example of a harmonic plane is based on the simplest harmonic metric 
$e^{2x}g_0$ and it can be shown that this harmonic plane is the Riemann surface for $\log z$ (see Section 7.1).  This Riemann surface is routinely constructed by gluing together an infinite number of Euclidean planes, cut along the negative x-axis.  The Euclidean metric on each plane extends smoothly across the glued edges, providing an “extended” Euclidean metric on the surface. The best visualization for this Riemann surface –-  though rarely (if ever) observed  -- is an infinitely spiraling cone in $\R^3$, which has the advantage of faithfully depicting the zero curvature everywhere (see Figure 1a).   

Under uniformization, when a biholomorphism exists mapping the Euclidean plane to a Riemann surface, the map is used to push-forward the Euclidean metric to the Riemann surface, thereby ``homogenizing" the surface.  Instead, we can use this map to pull-back the (extended) Euclidean metric, which arose from the cut-and-paste construction of the Riemann surface, thereby preserving the natural geometry of the construction.

\item
If $\varphi$ is a harmonic function on the unit disk, then the associated harmonic metric yields a \textit{harmonic disk}.  According to the Riemann mapping theorem, any proper simply-connected domain $\Omega$ in the Euclidean plane can be mapped biholomorphically to the unit disk.  If the Euclidean metric on $\Omega$ is pushed to the disk by the biholomorphism, then the resulting metric is a harmonic metric and the biholomorphism is an isometry.  Clearly, this harmonic metric is now a ``subset metric" of the Euclidean plane, in the sense that the harmonic disk can be isometrically embedded in the Euclidean plane.  Hence, there is a direct correspondence between proper simply-connected domains and harmonic “subset” metrics on the unit disk (up to pre-composition of the associated Riemann maps with rigid motions).  Importantly, there are harmonic metrics on the unit disk that are not subset metrics.



\end{enumerate}

Harmonic metrics are hiding in plain sight.
\section{Introduction}

\textit{Conventions: } We assume that all metrics $g$ are Riemannian and $C^\infty$ smooth unless otherwise noted.  Additionally, all surfaces are positively oriented, all mappings are orientation preserving, and a \textit{region} or \textit{domain} in the plane is a non-empty, open, connected subset. The notation $(\R^2,g)$ and $(D,g)$ represents the plane and the open disk, respectively, with metric $g.$ To ease exposition, we will use the function $F$ to represent a diffeomorphism and $f$ to represent a biholomorphism; of course, a diffeomorphism can be a biholomorphism.

\bigskip
Consider the (standard) Euclidean plane
$(\R^2,g_0)$.  By restricting the Euclidean metric to open proper subsets in the plane, one obtains a ready supply of -- seemingly inconsequential -- incomplete metrics. 
These are examples of what we call \textit{subset metrics}, 
which we define to be flat metrics on any open set that
can be realized as a proper subset of the Euclidean plane by 
an isometric embedding.  

Moreover, subset metrics can be constructed on the entire plane in (at least) two different ways.  If $F$ is a diffeomorphism from $\R^2$ to a proper open subset of $\R^2$, then the pull-back $F^*(g_0)$ of the Euclidean metric $g_0$  by $F$ gives $\R^2$ a flat metric for which $F$ is an isometry.  This realizes the pull-back metric as a subset metric.
Alternatively, subset metrics on $\R^2$ can be constructed as product metrics on $\R \times \R$; for example, consider the metric $e^{-2x^2} dx^2 + e^{-2y^2} dy^2$. This product metric can be realized as a subset of the Euclidean plane by an isometric embedding to an open square with sides of length $\sqrt{\pi}$.  However, the plane admits interesting flat metrics which are both natural and incomplete, but which are decidedly not subset or product metrics; the harmonic metrics described in Section 1 are examples.

The main purpose of this paper is to provide a classification of smooth incomplete flat Riemannian metrics on the plane (and the disk).  Sections 7 and 8 provide a diverse collection of examples and applications of harmonic planes and disks.

\section{Preliminaries}
It is well-known that a complete flat plane $(\R^2, g)$ is isometric to the Euclidean plane $(\R^2, g_0)$. Specifically, there is a diffeomorphism $F$ of the plane for which $g=F^*(g_0)$ is the pull-back metric

\begin{eqnarray*}
F:(\R^2,g) \rightarrow (\R^2,g_0)
\end{eqnarray*}
is an isometry.  This follows from the proof of the theorem of Cartan-Hadamard, where the exponential map $exp:  T_p(\R^2,g) \rightarrow (\R^2, g)$
is a diffeomorphism at any point $p\in \R^2.$  And by another theorem of Cartan on mappings that preserve curvature (see \cite{D}, Theorem 2.1, Chapter 8), it follows that this exponential map is an isometry. Using a linear isometry $L: T_p(\R^2,g) \rightarrow (\R^2, g_0)$, we set $F= L \circ exp^{-1}$ to conclude that, up to isometry, the Euclidean plane is the only complete flat plane.  If $F$ is also holomorphic, then it must be a complex affine transformation and 
$g=cg_0$, for some constant $c$. In the incomplete case, we have the following lemma (cf. \cite{BH}, Proposition 2.2).

\begin{lemma} \label{incomplete} 
   Every harmonic plane is incomplete.  
 \end{lemma}
 
\noindent
\textit{Proof:} Consider the harmonic plane $(\R^2, g_\varphi)$, with $g_\varphi = e^{2\varphi}g_0$ and 
$\varphi$ is a non-constant harmonic function.  Assume further, for a contradiction, that $(\R^2, g_\varphi)$ is complete. From the above discussion, there is a diffeomorphic isometry $F: (\R^2, g_\varphi)\rightarrow  (\R^2, g_0) $.   Since Euclidean angles are preserved, $F$ is holomorphic, hence entire (and affine), 
so $e^{2\varphi}= c$, for some constant c.  This is a contradiction. \qed

\begin{lemma} \label{another}
   No harmonic metric on the plane is a subset metric of the Euclidean plane.
 \end{lemma}
 
 \noindent
 \textit{Proof:}  Assume the contrary, so let  $F:  (\R^2,g_\psi) \rightarrow (S\subsetneq \R^2, g_0)$ be a (diffeomorphic) isometric embedding onto $S$, hence entire.  Picard's little theorem states that $S$ is either $\R^2$ or $\R^2-\{point\}$. It can't be the former since $S$ is a proper subset.  And it can't be the latter since $S$ is simply connected.  \qed
 
\bigskip
\noindent
This result can be immediately extended, as follows.

 \begin{lemma} \label{another2}
   No harmonic metric on the plane is a subset metric of another harmonic plane.
 \end{lemma}
 
 \noindent
 \textit{Proof:} The proof of Lemma~\ref{another} remains valid replacing $g_0$ with any harmonic metric. \qed

\section{The Classification Theorem(s)}



\begin{theorem}\label{thm:main} Every incomplete Riemannian flat plane $(\R^2, \widehat{g} )$ has the form $\widehat{g}=F^*(g)$ for exactly one of the following cases:

\bigskip

Case 1.  The map $F:  \R^2 \rightarrow (\R^2, g)$ is a diffeomorphism and $g$ is a harmonic metric,

\bigskip
Case 2. The map $F:  \R^2 \rightarrow (D, g)$ is a diffeomorphism and $g$ is a harmonic metric.


\bigskip
\noindent
In particular, the isometry class of $(\R^2, \widehat{g})$ contains a harmonic plane or a harmonic disk (but not both).  Furthermore,
for a harmonic metric $g$ on $\R^2$ or $D$, the isometry class $[g]$ in the space of harmonic metrics on $\R^2$ or $D$ is given by:

\bigskip
Case 1 \textup(continued\textup). The isometry class $[g]= \{f^*(g)~ |~f:  \R^2 \rightarrow (\R^2,g)$ is a complex affine mapping $z\rightarrow az+b         ,~for~ a,b \in \C, ~a\neq 0\}, $                   

\bigskip

Case 2 \textup(continued\textup). The isometry class $[g] = \{ f^*(g)~ |~f:  D \rightarrow (D,g)$ is a M{\"o}bius automorphism  $z\rightarrow e^{i\theta}       \frac{z-a}{\bar{a}z-1},~for~ |a|<1  ~and~ 0\le \theta <2\pi        \}.$  The Euclidean metric restricted to $D$ is an element in one of these isometry classes
(cf. Section 5.3).
\end{theorem}

\begin{remark} The classification statement for incomplete Riemannian flat metrics on the unit disk $D$ carries over
\textit{mutatis mutandis} from Theorem 4.1 -- with the obvious changes to the first sentence of the Theorem,  namely, replacing ``flat plane" $(\R^2, \widehat{g})$ with ``flat disk" $(D, \widehat{g})$ and replacing the domain $\R^2$ for the map $F$ in Cases 1 and 2 with $D$, and in the second sentence, replacing $(\R^2, \widehat{g})$ with $(D, \widehat{g})$.

\end{remark}

\section{Remarks on the Classification Theorem}

\subsection{Remarks on Case 1}
In this case, $F$ is holomorphic if and only if ~$\widehat{g}$ is a harmonic metric. In particular, there is no isometry from a harmonic plane to a harmonic disk (per Liouville).

In contrast to Lemmas \ref{another} and \ref{another2}, there do exist (local) \textit{isometric immersions} from harmonic planes to the Euclidean plane (and other harmonic planes).  For example, consider the 
``exponential mapping" $F(x,y)=(e^x\cos{y}, e^x\sin{y})$.  A straightforward calculation shows that $F^*(g_0)=e^{2x}g_0$, 
so this mapping is an isometric immersion $F:  (\R^2, e^{2x}g_0)\rightarrow (\R^2,g_0)$ with image $\R^2-\{origin\}$. See also Section 7.1.

\subsection{Remarks on Case 2}
In this case, we find that the geometry of a harmonic disk is more nuanced than the geometry of a harmonic plane.  In Section 1, harmonic disks were discussed in the context of the Riemann mapping theorem, where Riemann maps were used to create ``subset metric" disks and to provide the natural isometric embedding into Euclidean space.  In Section 7.1, we provide an easy example of a harmonic disk that cannot be isometrically embedded in the Euclidean plane; however, it does isometrically \textit{embed} in a harmonic plane 
$(   \R^2, e^{2x}g_0   )$.







If a harmonic disk cannot be isometrically embedded in any harmonic plane (or the Euclidean plane), and cannot be isometrically embedded as a proper subset in a harmonic disk (or the Euclidean disk), then we call it \textit{exotic}.  The existence of exotic harmonic disks has not yet been established, but our 
preliminary research on harmonic metrics derived from the real and imaginary parts of lacunary functions suggests the likelihood of their existence (cf. Section 7.4). 
  
\subsection{Additional Remarks}
  
The Euclidean metric on the disk $(D, g_0)$ requires special attention since $g_0$ is not a harmonic metric and hence its use is not allowed in Case 2.   This is in contrast to  Remark \ref{only remark}, where $(D, g_0)$ is used to construct an example of a flat plane based on the pull-back metric of a diffeomorphism $F: \R^2 \rightarrow (D, g_0)$.
However, the isometry class of the Euclidean metric on the disk is allowed in the classification theorem, as follows.  

If $\sigma: D \rightarrow D$ is a M{\"o}bius automorphism that is not a pure rotation, then $\sigma^*(g_0) = e^{2\varphi}g_0 = g_1$ is a harmonic metric on  $D$ (because $\varphi$ is non-constant).  Hence, for the isometry classes, we have $[g_0] = [g_1]$ and so $g_1$ and $[g_1]$ are allowed in Case 2 and its continuation.  

Another metric worth mentioning is the product metric $(\R, e^{2\alpha}) \times (\R, e^{2\beta})$, where $\alpha, \beta: \R \rightarrow \R$ are smooth functions.  The product is a flat plane, which is isometric to either the  entire plane (if both factors are complete) or a proper subset of the Euclidean plane (if either factor is incomplete).  


\section{Proof of the Classification Theorem}

\textit{Proof of Theorem \ref{thm:main}:}  Minding's theorem\footnote{Minding's namesake theorem of 1839 established that all surfaces having the same constant curvature must be locally isometric.  To the best of our knowledge, all published proofs of it, inclusive of Minding's original argument,
are existential in nature.   Unfortunately, Minding's paper is only available in the original German \textbf{\cite{Minding}} or in a Russian translation.  See \cite{CW0} for a constructive proof of Minding's theorem in the flat two-dimensional case which makes use of harmonic metrics.} guarantees the existence of a neighborhood around any point in $(\R^2, \widehat{g})$ which maps isometrically to an open set in the Euclidean plane.  
Consider an atlas composed of these ``Minding Maps".
The chart maps are isometries and these maps can be chosen so their images are open disks in the Euclidean plane.  The transition functions are necessarily restrictions of rigid motions.  Hence this ``Minding" atlas is real analytic, as is the flat Riemannian metric, and we denote the associated real-analytic flat surface by $(\R^2, \widehat{g})_\omega$. There is a natural inclusion 

$$
h_1: (\R^2, \widehat{g})_\omega \rightarrow (\R^2, \widehat{g}),
$$
which is both a $(C^\infty)$ diffeomorphism and an isometry.

Define an almost complex structure $J$ on $(\R^2, \widehat{g})_\omega$ as a length-preserving rotation by $90^{\circ}$, consistent with the orientation, and note that $J$ is real-analytic with vanishing Nijenhuis tensor (always true on oriented surfaces). In this real-analytic setting, the integration of $J$ produces a complex structure, compatible with $J$ and the flat metric $\widehat{g}$, thereby yielding a flat Riemann surface $(M, \widehat{g}  )$. This follows from the paper by Newlander and Nierenberg [13] (cf. the paper's first section that references the work of Frobenius, Libermann, Eckmann, and Frölicher on real-analytic manifolds) or it can be shown directly by complexifying the ``Minding" atlas and noting that the rigid motion transition functions are holomorphic.  Hence, there is a natural inclusion 

$$
h_2: (\R^2, \widehat{g})_\omega \rightarrow (M, \widehat{g}),
$$
which is both a $(C^\omega)$ diffeomorphism and an isometry.   

\noindent
\begin{remark}\label{only remark}
 One might infer from the above that 
$M$ is $\R^2$, but it is possible for $M$ to be $D$,    
as the following example illustrates.  Let $F:  \R^2 \rightarrow (D, g_0)$ be a smooth diffeomorphism and, 
using the pull-back metric, $(\R^2, F^*(g_0))$ becomes a flat surface, with $F$ an isometry.
The transition to a real-analytic structure and then to a complex-analytic structure transforms $(\R^2, F^*(g_0))$
into the Riemann surface $(D, g_0)$. 
\end{remark}

\noindent
By uniformization, there is a biholomorphishm $H$ between $\R^2$ (or $D$) and $M$. We leverage this biholomorphism to push the Euclidean $(x,y)$ coordinate system to a global coordinate system on $(M, \widehat{g})$, where the metric tensor 
$\widehat{g}$ has components:

\[
\widehat{g}_{12}=\widehat{g}_{21}=\widehat{g}(H_x,H_y)
\]
\[
\widehat{g}_{11}= \widehat{g}(H_x,H_x)~~and~~ \widehat{g}_{22}= \widehat{g}(H_y,H_y),
\]

\noindent
where the $x$ and $y$ subscripts are partial derivatives.

The compatibility of $J$ with the complex structure and metric on $(M, \widehat{g})$ means that $\widehat{g}_{12}=0$ and $\widehat{g}_{11}=\widehat{g}_{22}$.  Hence, the pull-back metric is $H^*(\widehat{g})=e^{2\varphi} g_0$, for the smooth real-valued function on $\R^2$ (or $D$) defined by $\varphi=\frac{1}{2}\log{(\widehat{g}_{11})}$.  Since $\widehat{g}$ is flat and incomplete, and $H$ is an isometry (in the pull-back metric), it follows that $\varphi$ is harmonic and not constant, i.e., $e^{2\varphi} g_0$ is a harmonic metric.  

Finally, the composition (for either $\R^2$ or $D$):

$$
G=h_1\circ h_2^{-1}\circ H: (\R^2 ~or~ D, e^{2\varphi}g_0) \rightarrow  (\R^2, \widehat{g})
$$

\noindent 
is a diffeomorphism and an isometry, and $G^{-1}$ is the required diffeomorphism $F$ in the theorem statement.  

We have shown that the isometry class of a flat plane $(\R^2,\widehat{g})$ includes a harmonic plane or harmonic disk.   Both cannot be included for if $F_1$ and $F_2$ are isometries from $(\R^2,\widehat{g})$ to $(\R^2, g_{\varphi_1})$ and $(D,g_{\varphi_2})$ respectively, the composition 

$$
F_2 \circ F_1^{-1}:  (\R^2, g_{\varphi_1}) \rightarrow (D, g_{\varphi_2})
$$
is an isometry, hence a conformal mapping -- with respect to Euclidean angles.  By complexifying $\R^2$ and $D$, the composition becomes a bounded non-constant holomorphic mapping -- a contradiction.

To complete the classification, we characterize the isometry classes of harmonic planes and harmonic disks as follows.  Let 
$f:  (\R^2, g_\varphi) \rightarrow (\R^2, g_\psi)$ be an isometric diffeomorphism between harmonic planes.  Then $g_\varphi= f^*(g_\psi)$ and, by complexifying, $f$ becomes a one-to-one entire function, hence a complex affine transformation. Similarly for harmonic disks, if $f:  (D,g_\varphi)\rightarrow (D, g_\psi)$ is an isometric diffeomorphism, then $f$ becomes a M{\"o}bius transformation. 
This completes the proof of the classification theorem.
\qed
\section{Examples and Applications}

This section presents examples and applications of the two primary constructions for harmonic metrics.  The two constructions are:

\begin{itemize}
    \item The pull-back harmonic metric $f^*(g_0)$ from a holomorphic function $f:\Omega \rightarrow \C$ with non-vanishing derivative (ensuring the metric is well-defined), and
    \item The harmonic metric formed directly from a harmonic function.
\end{itemize}

\noindent
On simply-connected domains $\Omega \subseteq \R^2$ , these two constructions are equivalent in the following sense.  Any harmonic metric of the form $e^{2\varphi}g_0$ can be realized as the pull-back metric $f^*(g_0)$, where $f:\Omega \rightarrow \C$ is defined by
 
\begin{eqnarray} \label{integral}
 f(z) = \int_{z_0}^z e^h dw,
\end{eqnarray}

 \noindent
 for any $z_0\in\Omega$ and a holomorphic $h:\Omega \rightarrow \C$
 with $Re(h)=\varphi$.  A straightforward calculation shows that
$f^*(g_0) =  |f^\prime(z)|^2g_0 = e^{2\varphi}g_0.$  

\bigskip
In Section 7.1, we consider the special case $h(z) = z$, and show the equivalence of the Riemann surface $\mathscr{L}$
for $\log{z}$ and the harmonic plane $(\R^2, e^{2x}g_0)$.  

In Section 7.2, we make direct use of the real and imaginary parts of the powers $z^n$, for $n > 1$.  Most of this section focuses on the construction of the Riemann surface identified with the harmonic plane having the harmonic metric derived from $Re(z^2) = x^2 - y^2$. 

In Section 7.3, we offer a variation of the usual cut-and-paste construction of Riemann surfaces and apply it to $\mathscr{L}$  discussed in Section 7.1.   Specifically, we replace the Euclidean plane building blocks (i.e., sheets of the Riemann surface) with harmonic planes and use this process to explore harmonic metrics with nested exponentials, such as $\varphi = Re(e^{e^z} +e^z +z).$

In Section 7.4, we consider lacunary functions on the unit disk, and leverage the harmonic disks associated with the real and imaginary parts to analyze the geometry near the unit circle boundary.  We focus on two lacunary functions and determine the length and curvature of selected radial line segments. 

\bigskip
\noindent
\subsection{The Simplest Harmonic Metric}
 Consider a biholomorphism $f: \mathscr{R}\rightarrow (\R^2,g_0)$ from a Riemann surface to the Euclidean plane.  The standard uniformization process uses the pull-back metric to install a complete flat metric on $\mathscr{R}$ and, in this manner, imposes the geometry of the Euclidean plane on all such Riemann surfaces.  As we will show, it is more geometrically intriguing to reverse this process and push-forward the (extended) euclidean metric $g_o$ on the Riemann surface $\mathscr{R}$ which is naturally inherited from each Euclidean plane ``building block", i.e., sheet, used in the cut-and-paste construction. In this way, the biholomorphism pushes the Euclidean metric to a harmonic metric on $\R^2$. In fact, all harmonic planes arise in a similar way,  as seen in Section 6, where the (rigid motion) transition maps of a holomorphic atlas define how Euclidean disks can be glued together to construct the associated Riemann surface with an extended $g_0$.

 Consider the Riemann surface $\mathscr{L}$ for $\log{z}$, consisting of an infinite stack of Euclidean planes, each cut along the positive $x$-axis with edges of adjacent planes glued together in the usual manner.  The surface $\mathscr{L}$ has a global polar coordinate system $(r,\theta)$ along with the extended Euclidean metric $g_0=dr^2+r^2d\theta^2$. We will use the biholomorphism $\log : \mathscr{L}\rightarrow \R^2$ defined by $\log (r,\theta)=(\log{r},\theta)$, where the ``log" function is overloaded to represent both the real log and the complex log in polar coordinates.

It is easy to show that the push-forward metric by the log mapping is $e^{2x}g_0$ on $\R^2$.  To analyze the geometry, we can leverage the Euclidean geometry of $\mathscr{L}$ to easily construct geometric objects and then push them from $\mathscr{L}$ to $(\R^2,e^{2x}g_0)$ by the log isometry.  For example, the geodesic rays in $\mathscr{L}$, defined by $\theta=$ constant and $0<r<\infty$, are isometrically mapped to geodesics in $(\R^2,e^{2x}g_0)$ as horizontal lines with finite length in the direction $x\rightarrow -\infty$;  see Figure 1a.  In $\mathscr{L}$, the geodesic rays converge as $r\rightarrow 0$, so in $(\R^2,e^{2x}g_0)$, the horizontal lines also converge as $x\rightarrow -\infty$.  Moreover, the constant curvature spirals in $\mathscr{L}$, corresponding to $r=$ constant and $-\infty<\theta<\infty$, isometrically map to the vertical lines in $(\R^2,e^{2x}g_0)$ and inherit the constant curvature $\frac{1}{r}$ from the corresponding spiral of radius $r$.  For example, the $y$-axis is distance one from $x=-\infty$ and hence has constant curvature one, which means the $y$-axis is an infinite-length, non-intersecting ``unit circle" (see Figure 1a).  This is one reason why $e^{2x}g_0$ is not a subset metric; it cannot be isometrically embedded as a subset of the Euclidean plane, though it can be isometrically immersed.

This metric also supports the construction of a harmonic disk that is not isometric to a subset of the Euclidean plane.  Let $\Omega$ be a simply-connected domain that contains the segment of the $y$-axis from $(0, 0)$ to $(0, 3\pi)$.  Using a Riemann map from $D$ to $\Omega$, the pull-back of $e^{2x}g_0  $ produces a harmonic disk that contains a non-intersecting curve with constant curvature one and length 3$\pi$.  Such a curve cannot exist in the Euclidean plane.

Of course, the metric $e^{2x}g_0$ is sufficiently simple for direct analysis, e.g., in $(x,y)$ coordinates, half of the Christoffel symbols are 0 and and the other half are $\pm 1$, so the geodesic equations are simple.  A visualization of $(\R^2, e^{2x}g_0)$ has already been described as the infinitely-sheeted Riemann surface for $\log{z}$, but a simpler visualization is an infinitely-spiraling cone (see Figure 1b).  Placing the cone's vertex at the origin in $\R^3$, we see that the intersection with the unit sphere is the ``unit circle".

$$
\includegraphics[height=3.2in]{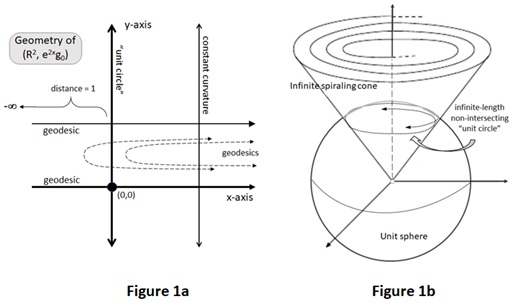}       
$$

\begin{footnotesize}
\noindent
\textbf{Figure 1:}  On the harmonic plane $(\R^2, e^{2x}g_0$), the behaviors of horizontal and vertical lines are labeled in Figure 1a.  The infinite spiraling cone in Figure 1b represents an isometric embedding of  $(\R^2, e^{2x}g_0)$ in $\R^3$, where the cone vertex corresponds to the point $x=-\infty$ in Figure 1a.
\end{footnotesize}

\bigskip
We end this example with a curious observation.  Consider the different locations of $e^{2x}$ in the following three metrics:  $e^{2x}(dx^2+dy^2)$, $e^{2x}dx^2 +dy^2$, $dx^2+e^{2x}dy^2$.  The first metric is the topic of this section. The second is a product metric and isometric to a Euclidean half-plane.  The third is a complete metric with constant curvature -1 and isometric to the upper half-plane with the Poincar{\' e} metric. 

\subsection{The Second Simplest Harmonic Metric}

In the previous section, we had a uniformization map from the Riemann surface for $\log{z}$ to the complex plane, and this map became an isometry when the (extended) Euclidean metric was pushed forward to construct the harmonic metric $e^{2x}g_0$.  This allowed the geometry of the harmonic metric to be understood in terms of the (Euclidean) geometry of the Riemann surface and the isometry map.  Without an explicit uniformization map or, equivalently, without an explicit solution to the integral (\ref{integral}) for a given harmonic function $\varphi$,  a push-forward metric $e^{2\varphi}g_0$ must be analyzed directly.

In this section, we consider the harmonic function Re$(z^2)=x^2-y^2$ and show that the geometry of $(\R^2,e^{2(x^2-y^2)})$ is surprisingly non-trivial. For simplicity, we only describe the behavior of radial lines $y=cx$ emanating from the origin;  see Figure 2, where we also provide a visualization in $\R^3$.  Figure 2a shows the first quadrant divided into two $45^\circ$ sectors.  The geometry of the diagonal ray separating these two sectors is an Euler spiral in the harmonic metric, i.e., its curvature is proportional to its arclength.  The rays corresponding 
to the positive $x$ and $y$ axes are geodesic rays, one with infinite length and the other with finite length 
$\frac{1}{2}\sqrt{\pi}$. The geometry of other first-quadrant rays is noted in the figure and this geometry is duplicated in the other three quadrants.  Figure 2b represents an isometric embedding of this harmonic plane in $\R^3$.  An Euler spiral (not shown) winds around each of the four infinite-spiraling cones.  All four cones are shown with initial up/down wrapping along the $y$-axis branch cut.  Each of the four Euler spirals converges to the cone vertex points of the form $(0, \pm \frac{1}{2}\sqrt{\pi}).$

\[
\includegraphics[height=3.2in]{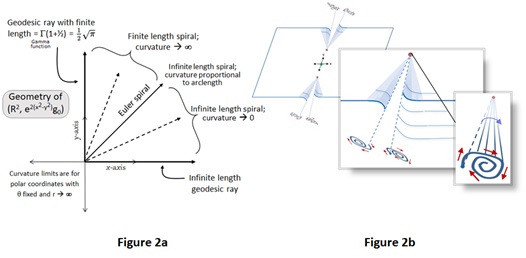}       
\]


\begin{footnotesize}
\noindent
\textbf{Figure 2:}
In Figure 2a, the behavior of rays in the first quadrant is shown for the specified harmonic metric.  Figure 2b shows a portion of the harmonic plane isometrically embedded in $\R^3$.  The two finite length geodesic rays (comprising the entire $y$-axis in Figure 2a) are correctly represented with finite length in Figure 2b, ending at the points $(0, \pm \frac{\sqrt{\pi}}{2})$, which are the start of two branch cuts from which four infinite spiraling cones emerge.
\end{footnotesize}

\bigskip
\bigskip
Finally, the behavior of the harmonic metrics based on Re$(z^n)$, for $n>2$, is similar.  There are $n$ finite-length geodesic rays, uniformly spaced around the origin, with length $\Gamma (1+\frac{1}{n})$, where $\Gamma$ is the Gamma function.  Between two adjacent finite-length geodesic rays there is an infinite length geodesic ray.  Between adjacent geodesic rays, there is a generalized Euler spiral; that is, one for which the curvature is directly proportional to a power of the arc length -- here the power is $n$.  The spirals between a geodesic ray and and the generalized Euler spiral follow the same pattern of behavior as described above for $n=2$. The isometric embedding in $\R^3$ has $n$ vertex points, and each of these points is the beginning of a branch cut which is a continuation of the direction of the finite-length geodesic.
The harmonic planes derived from Re$(z^n)$ and  Im$(z^n)$ are isometric by the rotation that sends the real part to the imaginary part, i.e., $\frac{\pi}{2n}$.
\subsection{Cut-and-Paste with Harmonic Planes}

In this section, we re-purpose harmonic planes to be the building blocks (replacing Euclidean planes) in the cut-and-paste construction of Riemann surfaces. We begin with a discussion of log-Riemann surfaces, as analyzed by Biswas and Perez-Marco in \cite{BP2}.  
In Section I.1.1 of their paper, these surfaces are defined by a cut-and-paste process using copies of the complex plane, where each plane/sheet can have one or more cuts, subject to various cutting and gluing conditions.  For example, all cuts are along straight rays, the set of rays are disjoint, and the set of ray vertices (i.e., ramification points) form a discrete set in each plane.  Furthermore, the number of cuts and the orientation of the cuts can differ between planes.  Importantly, all edges must have the same coordinate parameters prior to gluing so, for example, a cut edge along the positive x-axis cannot be glued to a cut edge along the positive y-axis.  This condition is key to replacing the Euclidean plane with a harmonic plane.

In Theorem II.5.3.1 of \cite{BP2}, it is shown that if there are a finite number of ramification points, each with infinite order, i.e., an infinite number of glued edges connect to each ramification point, then the Riemann surface is biholomorphic to $\C$. and the biholomorphism can be computed explicitly.  Other examples (see \cite{BP}) are available where different conditions can be imposed on the cut-and-paste construction to yield Riemann surfaces and explicit biholomorphisms to $\C$.  These biholomorphisms can then be used to push the extended Euclidean metric on the Riemann surfaces to produce harmonic metrics on $\C$.

So, for log-Riemann surfaces with biholomorphisms to $\C$, we can replace the Euclidean plane with any harmonic plane in the cut-and-paste process to create a ``harmonic" log-Riemann surface.  We provide an example of this process below.   It is worth noting that in  Section II.7 of \cite{BP}, the authors discuss some limitations of their methods when, for example, there are an infinite number of ramification points.   Perhaps some of these limitations could be addressed by leveraging harmonic planes as building blocks.  Since a harmonic plane already has one or more ramification points, the cut-and-paste construction for a log-Riemann surface will naturally yield a harmonic log-Riemann surface with an infinite number of ramification points.  The explicit biholomorphism for the initial log-Riemann surface remains valid for uniformization and it pushes the (extended) harmonic metric on the harmonic log-Riemann surface to a harmonic metric on $\C$.  

Our example begins with the usual $\log{z}$ Riemann surface $\mathscr{L}$, and replaces the Euclidean plane with the harmonic plane $(\R^2, e^{2x}g_0)$ which, as shown in Section 7.1, also represents $\mathscr{L}$.  The biholomorphism $\log: \mathscr{L}\rightarrow \R^2$ pushes forward the extended metric $e^{2x}g_0$ on $\mathscr{L}$ to $\R^2$, and a straight-forward calculation shows that the resulting harmonic metric is $e^{(e^x\cos{y}+x)}g_0$, with associated harmonic function Re($e^z +z$).

From Section 7.1 (and Figure 1a), we know the geometry of the building block $(\R^2, e^{2x}g_0)$.  And since the log function is an explicitly-defined isometry, we can push geometric objects in $(\R^2, e^{2x}g_0)$ to $(\R^2, e^{(e^x\cos{y}+x)}g_0)$.  This process is simple in concept, but the details are challenging and we cannot offer an isometric embedding in $\R^3$.

Furthermore, the process can be iterated, using the resulting harmonic plane at each stage as the building block for the next stage.  It is straight-forward to show that the second iteration of this example produces a harmonic metric on $\R^2$ with associated harmonic function Re($e^{e^z}+e^z +z)$.  The pattern is now apparent.
\subsection{Lacunary Functions}
The existence of analytic functions with natural boundaries -- that is, functions that cannot be extended analytically at any point on the circle of convergence -- was first discovered by Weierstrass and Kronecker in the 1860s.  Research into these \textit{lacunary functions} flourished through the mid-twentieth century, and continues today as an active area of study, based primarily on methods that are analytic in nature.  

In this section, we suggest that harmonic disks may provide a geometric basis for the study of these functions.  For example, using the harmonic disks associated with a lacunary function 
(via the real and imaginary parts), the behavior of various curves in these disks (e.g., radial lines, geodesics) offers insight into the geometry of lacunary functions near the $S^1$ boundary.

Below, we select two lacunary functions $f_1(z)$ and $f_2(z)$, and
describe the geometry of certain radial segments in the associated harmonic disks. 


\begin{eqnarray*}\label{Functions}
f_1(z)=\sum_1^\infty z^{2^n} ~ \text{and} ~ f_2(z)=\sum_1^\infty z^{n!}
\end{eqnarray*}


\noindent
Consider the two sets of angles $A$ and $B$:
 $$
 A=\{\theta\in[0,2\pi]~|~\frac{\theta}{2\pi}=\frac{p}{q},~ \textrm{with}~ p\in \N ~\textrm{and}~ q\in 2^s, ~\textrm{for}~ s\in \N\}, ~\textrm{and}~ B=\{\theta\in[0,2\pi]~|~\frac{\theta}{2\pi}\in \Q\}.
 $$
\noindent
For $f_1(z)$, it is a straightforward calculation (using equation 1.5 in \cite{Hu} for the curvature calculations) to establish the following geometry of radial line segments with angle $\theta \in A$ and radial parameter $r\in(0,1)$:

\begin{itemize}
\item Harmonic disk derived from Re$(f_1(z))$ 
\begin{itemize}
    \item Length of these radial segments is infinite
    \item Curvature of these radial segments approaches zero as $r\rightarrow 1$
\end{itemize}

\item Harmonic disk derived from Im$(f_2(z))$
\begin{itemize}
\item Length of these radial segments is finite
\item Curvature of these radial segments approaches $\infty$ as $r\rightarrow 1$

\end{itemize}

\end{itemize}

For $f_2(z)$, it is also a straightforward calculation to establish the following curvature properties of radial line segments with $\theta\in B$:

\begin{itemize}

\item  Harmonic disk derived from Re($f_2$), the curvature of radial segments approaches zero as $r\rightarrow 1$.

\item Harmonic disk derived from Im($f_2$), the curvature of radial segments approaches infinity as $r\rightarrow 1$.

\end{itemize}

Somewhat less straightforward is the fact that for the harmonic disks derived from Re($f_2$) and Im($f_2$), all radial segments have finite length. This latter fact follows from the curious convergence of the improper integral

$$I=\int_0^1\exp\left(\sum_{k=0}^{\infty}r^{k!}\right)~dr.$$



\section{Epilogue}
The impetus for this paper was the result of a calculation to determine the geometry of vertical and horizontal lines in the incomplete flat metric $e^{2x}g_0$.  As noted in Section 7.1, this led to the realization that the harmonic plane $(\R^2, e^{2x}g_0)$ is isometric to the Riemann surface for $\log{z}$ and that this surface could be isometrically embedded into $\R^3$ as a flat infinitely spiralling cone, in contrast to the more common renderings which are visibly non-flat.  We were surprised that such a nice observation, based on such a routine calculation, was absent from the literature.  Ultimately, we were led to the classification theorem of Section 4.  

The essential approach of this paper is to view the Riemannian geometry of flat surfaces through the lens of classical complex analysis.  It is particularly noteworthy that by avoiding uniformization, the richness of the geometry of incomplete flat metrics can  be exposed.  For example, in a recent note \cite{CW1}, we show that the well-known Four-Vertex Theorem is true for any flat plane. This is a non-trivial extension of the Four-Vertex Theorem since harmonic planes cannot be isometrically embedded into the Euclidean plane.  

These pleasing results, coupled with the importance of harmonic functions, suggest that the study of the geometry of harmonic metrics should not be ignored.  With the classification theorem in place, there are a number of questions which will be of ongoing interest.  We list a few.

\begin{enumerate}





\item The curve-shortening flows \cite{GH, G, H} start with a smooth Jordan curve and, during the flow, require that the intermediate curves remain smooth and without self-intersections.  Using the smooth Riemann Mapping Theorem \cite{BK},
the initial Jordan curve can be modelled as the $S^1$ boundary of a harmonic disk with harmonic metric that extends smoothly to the boundary.  During the flow, the intermediate curves can also be modelled as 
$S^1$ boundaries of harmonic disks, so the entire curve-shortening flow can be modelled as a flow of harmonic metrics (i.e., a flow of harmonic functions) on the closed unit disk. What are the defining equations for this flow of harmonic functions?

The authors were not able to construct the differential equation(s) for a harmonic flow that reproduced the classical curve shortening flow, but instead found and solved another differential equation for a harmonic flow (see equation (\ref{flow equation}) below).  Our approach relies on the observation that, for a given flow of harmonic metrics on the closed unit disk
$\bar{D}$, the associated flow of harmonic functions on $D$ can be characterized as a flow of Dirichlet initial condition on $S^1$.  

A harmonic flow on $D$ begins (at time 0) with a harmonic metric on $D$ obtained as the pull-back metric from a Riemann map $F$ of $D$ to the interior of the given smooth Jordan curve.  Since the Riemann map extends smoothly to the $S^1$ boundary of $D$, the  pull-back metric and associated harmonic function smoothly extend to $S^1$, thereby establishing the Dirichlet boundary conditions at time $0$.



Setting $\varphi(1,\theta , t)$ to be the flow of Dirichlet boundary conditions, with time parameter $t$ and $S^1$ parameters $(1,\theta)$, the flow differential equation is 
\begin{eqnarray}\label{flow equation}
\varphi_t(1,\theta, t) = -k(1,\theta, t)e^{\varphi (1,\theta , t)}+ 1,		
\end{eqnarray}
where $\varphi_t$ is the time derivative and $k(1, \theta, t)$ is the curvature of $S^1$ in the metric $g_{\varphi (r,\theta,t)}$ defined by the harmonic function $\varphi$, which has the specified boundary conditions.  The initial condition at $t = 0$ is written as $\varphi(1, \theta, 0) = \varphi (1, \theta )$ where, by abuse of notation, the second “$\varphi$” is derived from the Riemann map $F$.  

We can show that the flow solution to (\ref{flow equation}) converges to a circle and satisfies the condition that once the flow curve becomes convex, it remains convex. Furthermore, our flow solution includes explicit formulas for the length, area, and curvature of the flow curve for any time $t$.


\item
From the Nash Embedding Theorem, it is known that all harmonic planes and harmonic disks can be isometrically embedded into some Euclidean space.  What is the relationship between the harmonic functions and the dimension of the Euclidean space?  For harmonic disks being isometrically embedded into the Euclidean plane, this question is 
equivalent to asking for a characterization of harmonic metrics that are Euclidean subset metrics, and hence are associated to a Riemann map.  A related question asks for a characterization when the harmonic disk is convex.



\item Manifolds of Riemannian metrics have been studied by various researchers (e.g., \cite{MM}).  These infinite-dimensional manifolds have natural metrics with geodesics, Jacobi fields, etc.  For the manifold of harmonic metrics on the plane or unit disk, per our classification, what is the behavior of geodesics, Jacobi fields, and other geometric objects in terms of the underlying harmonic functions? For example, given a curve-shortening flow - along with its representation as an arc in the space of harmonic metrics on the disk - what is the geometry of this arc?

\end{enumerate}






\end{document}